# Optimal Scheduling of Distributed Energy Resources Considering Volt-VAr Controller of PV Smart Inverters

Zahra Soltani
*School of Electrical, Computer, and Energy Engineering*
*Arizona State University*
Tempe, Arizona

Shanshan Ma
*School of Electrical, Computer, and Energy Engineering*
*Arizona State University*
Tempe, Arizona

Mohammad Ghaljehei
*School of Electrical, Computer, and Energy Engineering*
*Arizona State University*
Tempe, Arizona

Mojdeh Khorsand
*School of Electrical, Computer, and Energy Engineering*
*Arizona State University*
Tempe, Arizona

*Abstract*—This paper proposes an operational scheduling model of distributed energy resources (DERs) and PV smart inverters with Volt-VAr controller using an accurate AC optimal power flow (ACOPF) in an unbalanced distribution network. A mathematical mixed-integer model of local Volt-VAr droop controller of the distributed mixed-phase PV smart inverters is proposed based on the IEEE 1547-2018 standard and is incorporated in the unbalanced ACOPF, which enables effective utilization of the Volt-VAr controllers to not only alleviate voltage issues locally but also at the feeder level. The proposed model is tested on two actual snapshots of a distribution feeder in Arizona. Also, the proposed operational scheduling method considering the Volt-VAr droop controller of PV smart inverters is compared with a recent work in scheduling of the PV smart inverters. The results illustrate that the PV smart inverters dispatches obtained by the proposed model can be practically implemented by local controller of inverters.

## Nomenclature

*Sets and Indices*
$\mathcal{K}/l$ — Set/Index of lines
$\mathcal{P}/i,j$ — Set/Index for buses
$\phi/x$ — Set/Index of phases
$F1/F2$ — Set of PV units with/without Volt-VAr controller

*Parameters and Constants*
$y_l^{x,m}/Z_l^{x,m}$ — Shunt admittance/Impedance between phases $x$ and $m$ of line $l$
$R_l^{x,m}/X_l^{x,m}$ — Resistance/Reactance between phases $x$ and $m$ of line $l$
$L_{o,x}^T$ — No load loss of transformer $o$ at phase $x$
$D_{d,x}^P/D_{d,x}^Q$ — Active/Reactive load $d$ at phase $x$
$M$ — A large positive number
$I_l^{max,x}$ — Maximum current of line $l$ at phase $x$
$P_{f,x}^{PVV,av}$ — Available active power of PV $f$ with Volt-VAr at phase $x$
$Q_{f,x}^{PVV,max}$ — Maximum reactive power of PV $f$ with Volt-VAr at phase $x$
$S_{f,x}^{PVV}$ — Apparent power rating of PV $f$ with Volt-VAr at phase $x$
$P_{e,x}^{PV,av}$ — Available active power of PV $e$ without Volt-VAr at phase $x$
$\pi_{s,x}^B/\pi_{r,x}^G$ — Substation $s$/DER $r$ electricity price
$\pi_{f,x}^{PVV}/\pi_{e,x}^{PV}$ — PV with/without Volt-VAr $f/e$ electricity price

*Variables*
$I_l^{r,x}/I_l^{im,x}$ — Real/Imaginary part of current flow at phase $x$ of line $l$
$V_i^x$ — Voltage magnitude at bus $i$ and phase $x$
$V_i^{r,x}/V_i^{im,x}$ — Real/Imaginary part of voltage at bus $i$ and phase $x$
$I_i^{r,x}/I_i^{im,x}$ — Real/Imaginary part of current injection at bus $i$ and phase $x$
$P_{r,x}^G/Q_{r,x}^G$ — Active/Reactive power of DER $r$ at phase $x$
$P_{f,x}^{PVV}/Q_{f,x}^{PVV}$ — Active/Reactive power of PV $f$ with Volt-VAr at phase $x$
$P_{s,x}^B/Q_{s,x}^B$ — Active/Reactive power of substation $s$ at phase $x$
$Q_{k,x}^C$ — Reactive power of capacitor $k$ at phase $x$

## I. Introduction

The fast voltage fluctuation and violation issues caused by high penetration of distributed energy resources (DERs), especially residential solar photovoltaic (PV) generation, cannot be easily coped with by traditional Volt-VAr control (VVC) devices, such as fixed/switched capacitor banks, load-tap changes, and voltage regulators [1]. The sparse number, limited operation frequency, and different working time scales of those traditional VVC devices result in a slow control reaction to the fast and local voltage violations. Furthermore, installing traditional VVC devices is costly [2]. To mitigate the fast and local voltage fluctuations and violations in the distribution systems, the smart inverter-based generators have been allowed to participate in the distribution feeder voltage regulation in amendment IEEE 1547-2018 standard by providing sufficient active and reactive power support [3]. Therefore, PV smart inverter with local voltage management capability is one of the most cost-effective VVC devices to regulate voltage issues at both local and system levels and increase the hosting capacity of distribution feeders.

Many studies have investigated different Volt-VAr optimization (VVO) approaches to manage the distribution feeders' voltage and reactive power profiles by dispatching PV smart inverters' active and reactive power output [4]-[7]. The authors in [1], [4], and [5] develop several centralized optimal

This work was supported by the Department of Energy Advanced Research Projects Agency – Energy under OPEN 2018 Program Award DE-AR0001001.

power flow (OPF) approaches for distribution system operator (DSO) to dispatch PV smart inverters outputs considering the voltage limit constraints. A distributed OPF approach based on alternating direction method of multipliers is proposed in [6] to dispatch the reactive power output of VVC devices to facilitate voltage reduction in unbalanced three-phase distribution systems. Although the studies in [1], [4]–[6] optimally dispatch the available active and reactive power outputs of PV smart inverters at the system level, they do not consider the PV smart inverter's internal Q-V droop control characteristics as defined in IEEE 1547-2018 [3]. Therefore, the scheduled reactive power from the OPF solution for a specific PV smart inverter may not be generated since the smart inverter's reactive power output varies with its local voltage magnitude. To dispatch the reactive power of PV smart inverters accurately to mitigate the fast voltage violation, authors in [7] integrate the piecewise model of Volt-VAr droop function (Q-V curve) as constraints into a centralized balanced distribution system OPF. However, the distribution system OPF with Volt/VAr droop function in [7] does not model the unbalanced distribution system. As a result, the dispatched reactive power output of a single-phase PV smart inverter may be inaccurate since the voltages are considered balanced in the OPF. Also, no comparison with the prior methods in the distribution system OPF without Volt/VAr droop modeling is provided in [7] to illustrate the necessity of modeling Volt/VAr droop function.

This paper proposes a novel distribution system operational scheduling strategy to optimally schedule the DERs with the objective of minimizing the total operational cost while keeping the system voltage at the acceptable range. The proposed scheduling strategy is modeled as an unbalanced ACOPF based on current and voltage (IVACOPF) formulation and includes a detailed representation of Volt-VAr droop controller of mixed-phase distributed PV smart inverters. A set of novel mixed-integer linear constraints are developed to formulate the Q-V characteristics of Volt-VAr controllers of distributed PV smart inverters using the Big-M method to present the piecewise Volt-VAr droop function based on the IEEE 1547-2018 standard. The proposed model is tested on a real large primary distribution feeder in Arizona using actual data from two snapshots. The results illustrate the effectiveness of incorporating the distributed PV smart inverters' Volt-VAr droop control function in distribution system scheduling for achieving local as well as system-level voltage regulation in an unbalanced distribution system. Also, the performance of the proposed model is compared with a prior work [5], which models the PV smart inverters without considering their droop settings.

This paper is organized as follows: Section II describes the mathematical formulation of the proposed scheduling problem. Section III conducts the simulation, and Section IV summarizes the conclusion.

## II. PROBLEM FORMULATION

First, the formulation of optimal dispatchable resource scheduling using IVACOPF is presented. Then, the proposed model of integrating accurate Vol-VAr controller of the smart PV inverters in the resource scheduling problem is discussed.

### A. Optimal DERs Scheduling Using IVACOPF in an Unbalanced Distribution System

Distribution systems are three-phase unbalanced networks with high R to X ratio of distribution lines. In order to solve ACOPF in a distribution system, DistFlow model is widely used in the literature [8]. However, the DistFlow model disregards line losses and assumes that bus voltages are balanced. In this paper, a new ACOPF model based on current and voltage (IVACOPF) is used, which models all details of distribution systems including untransposed distribution lines with mutual impedances and mutual admittances as well as line losses [9]. The IVACOPF model for an unbalanced distribution system is in rectangular coordinate. In IVACOPF, the current in one of the phases of a three-phase distribution line from bus $i$ to $j$ bus is formulated as (1).

$$I_l^x = (Z_l^{x,x})^{-1}\left[V_i^x - V_j^x - \sum_{m\in\phi, m\neq x} Z_l^{x,m}(I_l^m) + \frac{1}{2}\sum_{m\in\phi} Z_l^{x,m}\left(\sum_{n\in\phi} y_l^{m,n} V_i^n\right)\right], \forall\, x \in \phi, l \in \mathcal{K} \quad (1)$$

It can be seen in (1) that the current in phase $x$ depends on voltages of phase $x$ and also the voltage and current of other phases because of mutual impedances and admittances. Constraint (1) in a rectangular coordinate is defined as (2)-(3). It is worth noting that (2)-(3) are linear constraints.

$$I_l^{r,x} = (R_l^{x,x})^{-1}\left[V_i^{r,x} - V_j^{r,x} - \sum_{m\in\phi, m\neq x} R_l^{x,m} I_l^{r,m} - \frac{1}{2}\sum_{m\in\phi} R_l^{x,m}\left(\sum_{n\in\phi} y_l^{m,n} V_i^{im,n}\right) + \sum_{m\in\phi} X_l^{x,m}\left(I_l^{im,m} - \frac{1}{2}\sum_{n\in\phi} y_l^{m,n} V_i^{r,n}\right)\right], \forall\, x \in \phi, l \in \mathcal{K} \quad (2)$$

$$I_l^{im,x} = (R_l^{x,x})^{-1}\left[V_i^{im,x} - V_j^{im,x} - \sum_{m\in\phi, m\neq x} R_l^{x,m} I_l^{im,m} + \frac{1}{2}\sum_{m\in\phi} R_l^{x,m}\left(\sum_{n\in\phi} y_l^{m,n} V_i^{r,n}\right) - \sum_{m\in\phi} X_l^{x,m}\left(I_l^{r,m} + \frac{1}{2}\sum_{n\in\phi} y_l^{m,n} V_i^{im,n}\right)\right], \forall\, x \in \phi, l \in \mathcal{K} \quad (3)$$

At each bus of a distribution network, the injected current at each phase is formulated using (4)-(5) in a rectangular coordinate.

$$I_i^{r,x} = \sum_{l\in\hbar(i)} I_l^{r,x}, \forall\, x \in \phi, i \in \mathcal{P} \quad (4)$$

$$I_i^{im,x} = \sum_{l\in\hbar(i)} I_l^{im,x}, \forall\, x \in \phi, i \in \mathcal{P} \quad (5)$$

where $\hbar(i)$ is set of connected lines to bus $i$. The active power balance constraint for each phase of a bus in an unbalanced distribution system is given in (6). The set of substations, DERs, PV units with Volt-VAr controller, PV units without Volt-VAr controller, load, and distribution transformers that are connected to bus $i$ are presented by $s(i), r(i), f1(i), f2(i), d(i)$, and $s(i)$, respectively. The output active power of the DERs including battery units and distributed generations (DGs), the active power of PV units with and without Volt-VAr controller, and no-load loss of distribution transformers are modeled in (6).

$$\sum_{\forall s\in s(i)} P_{s,x}^B + \sum_{\forall r\in r(i)} P_{r,x}^G + \sum_{\forall f\in f1(i)} P_{f,x}^{PVV} + \sum_{\forall e\in f2(i)} P_{e,x}^{PV} - \sum_{\forall d\in d(i)} D_{d,x}^P - \sum_{\forall o\in o(n)} L_{o,x}^T = P_{i,x} = V_i^{r,x} I_i^{r,x} + V_i^{im,x} I_i^{im,x}, \forall\, x \in \phi, i \in \mathcal{P} \quad (6)$$

Constraint (7) models the reactive power balance constraint in each phase of a bus including reactive output power of DERs, PV with Volt-VAr, and capacitors. The set of connected capacitors to bus $i$ is given by $k(i)$.

$$\sum_{\forall s \epsilon S(i)} Q_{s,x}^B + \sum_{\forall r \epsilon r(i)} Q_{r,x}^G + \sum_{\forall f \epsilon f1(i)} Q_{f,x}^{PVV} + \sum_{k \forall k(i)} Q_{k,x}^C - \sum_{\forall d \epsilon d(i)} D_{d,x}^Q = Q_{i,x} = V_i^{im,x} I_i^{r,x} - V_i^{r,x} I_i^{im,x}, \forall\, x \in \phi, i \in \mathcal{P} \quad (7)$$

Constraints (6)-(7) are nonconvex and nonlinear. Using first-order approximation of Taylor series, (6)-(7) are reformulated as linear constraints (8)-(9) around $(V_i^{r,x(t)}, V_i^{im,x(t)})$ and $(I_i^{r,x(t)}, I_i^{im,x(t)})$.

$$P_{i,x} = V_i^{r,x(t)} I_i^{r,x} + V_i^{im,x(t)} I_i^{im,x} + I_i^{r,x(t)} V_i^{r,x} + I_i^{im,x(t)} V_i^{im,x} - V_i^{r,x(t)} I_i^{r,x(t)} - V_i^{im,x(t)} I_i^{im,x(t)}, \forall\, x \in \phi, i \in \mathcal{P} \quad (8)$$

$$Q_{i,x} = V_i^{im,x(t)} I_i^{r,x} - V_i^{r,x(t)} I_i^{im,x} + I_i^{r,x(t)} V_i^{im,x} - I_i^{im,x(t)} V_i^{r,x} - V_i^{im,x(t)} I_i^{r,x(t)} + V_i^{r,x(t)} I_i^{im,x(t)}, \forall\, x \in \phi, i \in \mathcal{P} \quad (9)$$

To this end, the model is solved iteratively and the parameters of the Taylor series, i.e., $V_i^{r,x(t)}, V_i^{im,x(t)}, I_i^{r,x(t)}, I_i^{im,x(t)}$, are updated in each iteration based on results of the previous iteration. Since the IVACOPF is formulated based on the rectangular coordinate, the voltage magnitude in each phase of a bus is defined as nonlinear constraint (10).

$$V_i^x = \sqrt{V_i^{r,x2} + V_i^{im,x2}}, \forall\, x \in \phi, i \in \mathcal{P} \quad (10)$$

A linear approximation of (10) is developed in this paper utilizing Taylor series with ignoring orders greater than one as defined in (11).

$$V_i^x = \frac{V_i^{r,x(t)}}{\sqrt{V_i^{r,x(t)^2} + V_i^{im,x(t)^2}}} V_i^{r,x} + \frac{V_i^{im,x(t)}}{\sqrt{V_i^{r,x(t)^2} + V_i^{im,x(t)^2}}} V_i^{im,x}, \forall\, x \in \phi, i \in \mathcal{P} \quad (11)$$

The limit for voltage magnitude of each phase of a bus is given in (12).

$$0.95 \leq V_i^x \leq 1.05, \forall\, x \in \phi, i \in \mathcal{P} \quad (12)$$

The thermal line limit of each phase of a three-phase distribution line is defined as (13). The output of a capacitor bank units in the distribution system is restricted by (14).

$$I_l^{r,x2} + I_l^{im,x2} \leq I_l^{max,x2}, \forall\, x \in \phi, l \in \mathcal{K} \quad (13)$$

$$Q_{k,x}^C \leq Q_{k,x}^{C,max}, \forall\, x \in \phi, k \in \mathcal{C} \quad (14)$$

The PV units without Volt-VAr controller only generates active power and are modeled using (15). The active and reactive power output of a PV unit with Volt-VAr controller is limited based on its apparent power rating given in (16). Also, active and reactive power outputs of a PV unit with Volt-VAr controller is restricted by the available active power of the PV unit based on the solar radiation and the apparent power rating, respectively, as shown in (17)-(18).

$$P_{e,x}^{PV} = P_{e,x}^{PV,av}, \forall\, x \in \phi, e \in F2 \quad (15)$$

$$Q_{f,x}^{PVV^2} + P_{f,x}^{PVV^2} \leq S_{f,x}^{PVV^2}, \forall\, x \in \phi, f \in F1 \quad (16)$$

$$P_{f,x}^{PVV} \leq P_{f,x}^{PVV,av}, \forall\, x \in \phi, f \in F1 \quad (17)$$

$$-S_{f,x}^{PVV} \leq Q_{f,x}^{PVV} \leq S_{f,x}^{PVV}, \forall\, x \in \phi, f \in F1 \quad (18)$$

The objective function of the dispatchable resource scheduling in an unbalanced distribution system in a general form is presented in (19), which minimizes total distribution system operation cost including cost of purchasing power from the bulk grid, DERs, and PV units with Volt-VAr controller. In (19), $R$ and $S$ are set of all the DERs and the substations in the system, respectively. The objective function in (19) is subject to constraints given in (2)-(5), (8)-(9), and (11)-(18).

$$\min_{\substack{V_i^{r,x}, V_i^{im,x}, I_i^{r,x}, I_i^{im,x}, I_l^{r,x}, \\ I_l^{im,x}, P_{r,x}^G, P_{e,x}^{PV}, P_{f,x}^{PVV}, \\ Q_{r,x}^G, Q_{f,x}^{PVV}, Q_{k,x}^C, V_i^x}} \left\{ \sum_{x \in \phi} \begin{matrix} \sum_{\forall s \in S} \pi_{s,x}^B P_{s,x}^B + \sum_{\forall r \in R} \pi_{r,x}^G P_{r,x}^G \\ + \sum_{\forall f \in F1} \pi_{f,x}^{PVV} P_{f,x}^{PVV} \end{matrix} \right\} \quad (19)$$

### B. Proposed DERs Scheduling Model Considering Volt-VAr Droop Controller of Smart PV Inverters

High penetration of the DERs such as behind-the-meter PV units in the modern distribution system may lead to voltage issues. According to the IEEE 1547-2018 standard, the DERs can contribute to regulating the voltage in the distribution system by injecting/absorbing the reactive power. Smart inverters of the PV units can mitigate voltage fluctuations by providing the reactive power support, which can result in less active power curtailment of the PV units by maintaining the system voltage within the acceptable range. Most of the prior work models the PV smart inverter only using (16)-(18) without considering the Volt-VAr droop controller of the inverter [5]. However, not modeling the inverter's droop Volt-VAr controller precisely may lead to an inaccurate dispatch, which may not follow the local inverter setting based on IEEE 1547-2018 standard. According to the IEEE 1547-2018 standard, Volt-VAr controller of a PV smart inverter should function based on the Q-V curve of Fig. 1.

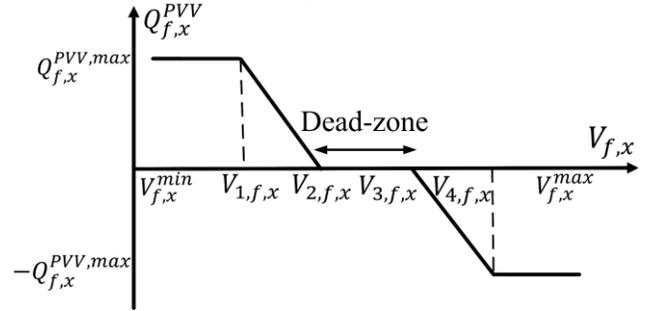

Fig. 1. Q-V characteristic of a PV smart inverter with Volt-VAr controller [3]

The voltage-reactive controller mode of PV units shown in Fig. 1 can be mathematically modeled as (20).

$$Q_{f,x}^{PVV} = \begin{cases} Q_{f,x}^{PVV,max} & V_{f,x}^{min} \leq V_{f,x} < V_{1,f,x} \\ \frac{Q_{f,x}^{PVV,max}}{V_{2,f,x}-V_{1,f,x}}(V_{2,f,x} - V_{f,x}) & V_{1,f,x} \leq V_{f,x} < V_{2,f,x} \\ 0 & V_{2,f,x} \leq V_{f,x} < V_{3,f,x} \\ \frac{Q_{f,x}^{PVV,max}}{V_{3,f,x}-V_{4,f,x}}(V_{f,x} - V_{3,f,x}) & V_{3,f,x} \leq V_{f,x} < V_{4,f,x} \\ -Q_{f,x}^{PVV,max} & V_{4,f,x} \leq V_{f,x} < V_{f,x}^{max} \end{cases} \quad (20)$$

where $V_{1,f,x}, V_{2,f,x}, V_{3,f,x}$, and $V_{4,f,x}$ are set points of the Q-V curve of a PV smart inverter. Equation (20) is a piecewise linear function. Based on voltage magnitude at a PV location (i.e., $V_{i,x}$), the controller of the PV smart inverter can operate in five different operating zone. This paper models (20) in an ACOPF to effectively account for the Q-V characteristics of Volt-VAr controllers while scheduling dispatchable resources of a distribution network. However, (20) includes "if-then" constraints, which cannot be simply modeled in the optimization framework. To this end, a mixed-integer linear programming formulation of (20) is proposed in this paper using the Big-M method in which each operation zones of Q-V curve of Fig. 1 is modeled using binary variables ( i.e., $z_{1,f,x}$, $z_{2,f,x}$, $z_{3,f,x}$, $z_{4,f,x}$, $z_{5,f,x}$). The proposed mixed-integer linear model of the Q-V curve of a PV smart invert is presented in (21)-(31).

$$V_{f,x} \leq V_{1,f,x} + Mz_{1,f,x} \quad (21)$$
$$-Mz_{1,f,x} \leq Q_{f,x}^{PVV} - Q_{f,x}^{PVV,max} \leq Mz_{1,f,x} \quad (22)$$
$$-Mz_{2,f,x} + V_{1,f,x} \leq V_{f,x} \leq V_{2,f,x} + Mz_{2,f,x} \quad (23)$$
$$-Mz_{2,f,x} \leq Q_{f,x}^{PVV} - \frac{Q_{f,x}^{PVV,max}}{V_{2,f,x} - V_{1,f,x}}(V_{2,f,x} - V_{f,x}) \leq Mz_{2,f,x} \quad (24)$$
$$-Mz_{3,f,x} + V_{2,f,x} \leq V_{f,x} \leq V_{3,f,x} + Mz_{3,f,x} \quad (25)$$
$$-Mz_{3,f,x} \leq Q_{f,x}^{PVV} \leq Mz_{3,f,x} \quad (26)$$
$$-Mz_{4,f,x} + V_{3,f,x} \leq V_{f,x} \leq V_{4,f,x} + Mz_{4,f,x} \quad (27)$$
$$-Mz_{4,f,x} \leq Q_{f,x}^{PVV} - \frac{Q_{f,x}^{PVV,max}}{V_{3,f,x} - V_{4,f,x}}(V_{f,x} - V_{3,f,x}) \leq Mz_{4,f,x} \quad (28)$$
$$-Mz_{5,f,x} + V_{4,f,x} \leq V_{f,x} \quad (29)$$
$$-Mz_{5,f,x} \leq Q_{f,x}^{PVV} + Q_{f,x}^{PVV,max} \leq Mz_{5,f,x} \quad (30)$$
$$z_{1,f,x} + z_{2,f,x} + z_{3,f,x} + z_{4,f,x} + z_{5,f,x} \leq 4 \quad (31)$$

Constraints (21)-(22) model the first operating zone of the voltage-reactive curve, where the PV unit injects maximum reactive power (i.e., $Q_{f,x}^{PVV,max}$) to the network. The injected reactive power of the PV unit at zone 2 is formulated using (23)-(24) as a linear function of local voltage magnitude at the PV location (i.e., $V_{f,x}$) with a positive slope. In the zone 3, where $V_{f,x}$ falls in the dead-zone, the reactive power of the PV unit is equal to zero as expressed in (25)-(26). If the local voltage at the PV location is larger than $V_{3,f,x}$, the PV unit will operate in inductive mode, i.e., zones 4 or 5 in Fig. 1, and will absorb the reactive power as modeled using (27)-(30). Constraint (31) forces the controller to operate in only one operating zone at each time instant. It should be noted that the PV unit operates in one of the five zones, if the binary variable associated with that zone is equals to zero. For instance, if $V_{f,x}$ is in the dead zone between $V_{2,f,x}$ and $V_{3,f,x}$, then $z_{3,f,x}$ will be obtained as zero according to (25) and as a result the reactive power of the PV unit will be equal to zero using (26).

The proposed detailed model of Volt-VAr droop controller of PV smart inverters, i.e., (21)-(31), is integrated with the dispatchable resource scheduling method presented in Section II. A. The objective function of the proposed dispatchable resource scheduling problem with accurate modeling of smart inverter of PV units with Volt-VAr controller is presented in (19), which is subject to constraints ((2)-(5), (8)-(9), (11)-(18), and also proposed constraints (21)-(31). The proposed DERs scheduling problem is a mixed-integer quadratically constrained program (MIQCP) model, which can be easily solved with commercial solver such as CPLEX and GUROBI.

## III. SIMULATION RESULTS

The proposed model is tested using actual system data of a primary distribution feeder of a utility in Arizona. The utility feeder has 2100 buses and 1790 lines [9]. There are 371 distribution transformers, 342 aggregated load buses, and 249 PV units while 77 PV units are equipped with Volt-VAr controllers. The test system represents actual system condition on March 15, 2019, 2 pm, with high PV penetration about 232% i.e., 3625 kW PV/1563 kW load) of the total feeder load, which contains overvoltage issues. Set points of PV units (i.e., $V_{1,f,x}$, $V_{2,f,x}$, $V_{3,f,x}$, and $V_{4,f,x}$) are considered to be 0.94, 0.98, 1.02, and 1.06 based on IEEE standard 1547-2018 [3]. $Q_{f,x}^{PVV,max}$ is defined based on the rating of each PV unit. The proposed model is simulated in JAVA and solved by GUROBI solver. Two case studies are considered (1) without Volt-VAr controllers of PV units (all 249 PV units without Volt-VAr controller) (2) with considering Volt-VAr controllers of PV units (77 PV units with Volt-VAr controller and 172 PV units without Volt-VAr controller). Table I shows the average error of the proposed linear approximation of the voltage magnitude (i.e., constraint (11)) among all buses and phases of the utility feeder in cases 1 and 2. It should be noted that the proposed IVACOPF model converges in two iterations in both cases. It can be seen in Table I that the proposed linear approximation of (10) is accurate in both cases with only two iterations.

Table I: Voltage magnitude approximation error

| Iteration | Case 1 (pu) | Case 2 (pu) |
|---|---|---|
| 1 | 1.28e-03 | 1.27e-03 |
| 2 | 2.10e-07 | 2.23e-07 |

The reactive power output of PV units and their corresponding local voltage magnitude in each of the 77 PV locations with and without considering Volt-VAr controller for these PV units ,i.e., cases 1 and 2, are shown in Fig. 2. It can be seen in Fig. 2 that without Volt-VAr controller, the local voltage magnitude in some of the PVs locations violates the allowable system limit ( i.e., 1.05 pu). However, with considering the proposed Volt-VAr optimization model of the PV smart inverters, all 77 PV units absorb reactive power and mitigate the local overvoltage issue. Voltage profile of all nodes of the utility feeder without and with considering Volt-VAr controllers of the PV units are shown in Fig. 3 and Fig. 4,

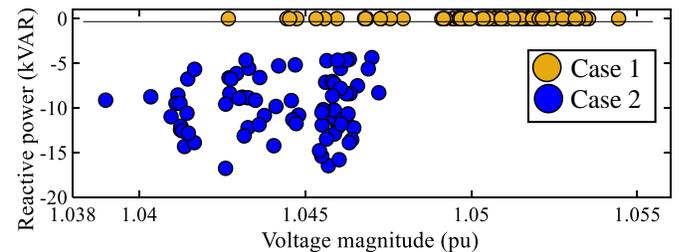

Fig. 2. Reactive power output of 77 PV units and their corresponding local voltage magnitude in cases 1 and 2.

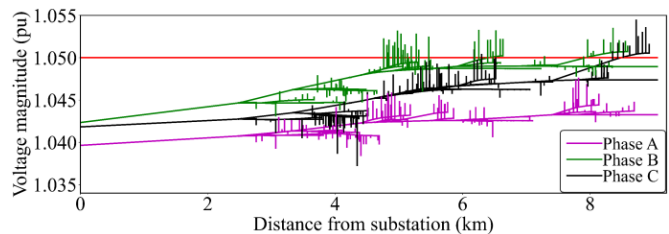

Fig. 3. Voltage of all nodes of the utility feeder in case 1 without Volt-VAr controllers.

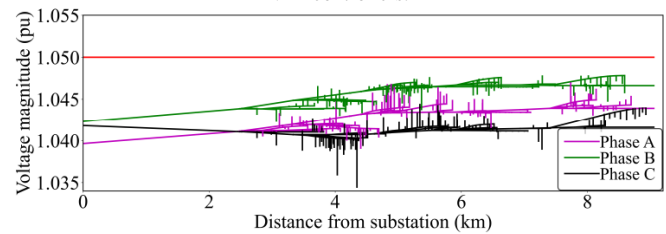

Fig. 4. Voltage of all nodes of the utility feeder in case 2 with considering Volt-VAr controllers for the 77 PV units.

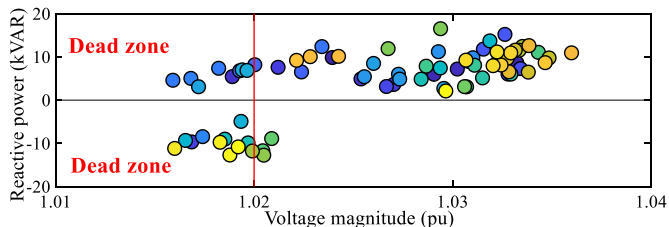

Fig. 5. Reactive power and voltage of PV units with Volt-VAr controller based on model of [5].

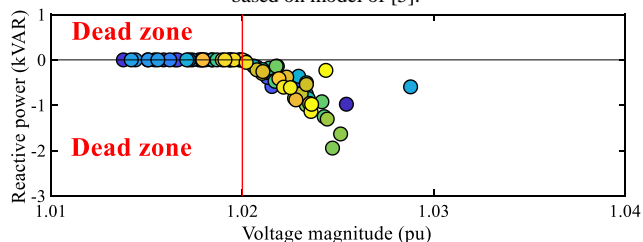

Fig. 6. Reactive power and voltage of PV units with Volt-VAr controller based on the proposed model

, respectively. By comparing Fig. 3 and Fig. 4, it can be seen that by considering Volt-VAr controller of the PV smart inverters, not only the voltage violation is resolved locally (shown in Fig. 2) but also the system-level overvoltage issue is resolved and the voltage of all nodes is within the allowable range due to the local reactive power support of the PV smart inverters.

*A. Comparison of the Proposed Method with Prior Work*

In this section, the performance of the proposed model is compared with prior work [5] in the modeling of the smart inverter of the PV units. In [5], PV smart inverters are modeled only using (16)-(18) without considering droop controller of PV smart inverters. In this regard, the proposed DERs scheduling model presented in section II.B with detailed representation of droop controller of PV smart inverters (i.e., constraints (16)-(18) and (21)-(31)) is compared with the DERs scheduling model in Section II.A, where smart inverter of PV units are modeled only by constraints (16)-(18) [5]. For the purpose of this comparison, a snapshot of the utility system on July 15, 2019, 5 pm with high load and low PV penetration is considered. Figure 5 shows the reactive power output of 77 PV units and the corresponding local voltage magnitude at each of these PV locations. It can be seen in Fig. 5 that the local voltage magnitudes in some of the PV locations is between 0.98 pu and 1.02 pu, which is the dead-zone of smart inverters' Volt-VAr controller; thus, their reactive power output must be equal to zero according to the IEEE 1547-2018 standards. For the rest of the PV units, the local voltage magnitude is larger than 1.02 pu, which is the zone 4 of the Q-V curve of droop controller of the PV smart inverters and their reactive power must be negative, i.e., absorption of reactive power, according to Fig. 1. However, since the method of [5] does not model details of the droop controller of the PV smart inverters based on the IEEE 1547-2018 standards, the solution of the DERs scheduling model for the PV smart inverters' dispatch is not accurate and practical in the local setting of the inverters. Figure 6 shows the reactive power output of the PV smart inverters with Volt-VAr controller based on the proposed model and their corresponding local voltage magnitude. It can be seen in the Fig. 6 that the proposed model schedules the PV smart inverters based on IEEE 1547 standards, which is a practical and precise solution for real-time implementation in the local inverter setting.

IV. CONCLUSION

In this paper, an operational scheduling framework is proposed for mitigating the voltage violations in the unbalanced distribution systems with DERs and PV smart inverters using an unbalanced distribution system ACOPF. The Volt-VAr droop controller of the mixed-phase distributed PV smart inverters is modeled mathematically using the Big-M method. The proposed Q-V characteristic of the PV units with Vol-VAr controller is integrated in the ACOPF model of an unbalanced distribution system to more accurately schedule these resources. A large distribution feeder in Arizona is utilized for testing the proposed model in two different actual snapshots of the system with high and low PV penetrations. The results illustrate the effectiveness of the proposed model for enhancing the system reliability and security by the accurate scheduling of the distributed PV smart inverters and enabling the local reactive power support for the local and the feeder-level voltage regulation. Also, a comparison case study with a recent work is conducted to further illustrates the performance of the proposed model.